\documentclass[12pt]{amsart}

\usepackage{color} 
\usepackage{graphicx} 
\usepackage{tabularx}
\usepackage{latexsym,amssymb,amsmath,amsfonts,amscd} 
\usepackage[all]{xy}
\usepackage{mathrsfs} 
\usepackage{tikz} 
\usetikzlibrary{shapes,patterns,arrows,decorations.pathreplacing}
\usepackage{epigraph,enumerate}
\usepackage{stmaryrd}
\usepackage{tikz-cd}
\usepackage{mathtools}
\usepackage{hyperref}

\allowdisplaybreaks[4]
\usepackage{latexsym, amssymb, amsmath}
\usepackage{enumerate}
\usepackage{flushend} 
\usepackage{mathrsfs} 
\usepackage{stfloats}
\def\custombibliography#1{
 \normalsize
\section*{\centering References}
 \list
 {[\arabic{enumi}]}{\settowidth\labelwidth{[#1]}\leftmargin\labelwidth
 \setlength{\itemsep}{.1em}
 \advance\leftmargin\labelsep
 \usecounter{enumi}}
 \def\newblock{\hskip .11em plus .33em minus -.07em}
 \sloppy
 \sfcode`\.=1000\relax}

\def\L2{{\cal L}_2}

\newcommand\bull{\vrule height .9ex width .8ex depth -.1ex } 

\newcommand\re{\rm I\! R}
\newcommand\cdcout[1]{} 

\newcommand{\rv}[1]{\boldsymbol{#1}} 
\newcommand{\RomanNumber}[1]{\uppercase\expandafter{\romannumeral #1}}
\newcommand{\romannumber}[1]{\lowercase\expandafter{\romannumeral #1}}

\DeclareMathAlphabet{\mathpzc}{OT1}{pzc}{m}{it}

\def\1{\rv 1} 

\usepackage{flushend} 
\usepackage{stfloats}
\usepackage{color} 
\usepackage{latexsym,amssymb,amsmath} 
\usepackage{mathrsfs,mathtools} 
\usepackage{multicol,enumerate}
\usepackage{subcaption}

\def\abs#1{\lvert #1 \rvert}

\def\allpoly{\mbox{$\re\langle X \rangle$}}

\def\allpolyx0degn{\mbox{$P_n$}}

\def\allwords{\mbox{$X^{\ast}$}}

\def\allseries{\mbox{$\re\langle\langle X \rangle\rangle$}}
\def\allseries#1{\mbox{$\re^{#1}\langle\langle X \rangle\rangle$}}

\def\allproperseries#1{\mbox{$\re_{p}^{#1}\, \langle\langle X \rangle\rangle$}}

\def\allseriesX1{\mbox{$\re [[ X_1 ]]$}}

\def\bull{\rule{0.08in}{0.08in}} 

\newcommand{\comment}[1]{} 

\def\eqref#1{(\ref{#1})} 

\def\lhook{\curvearrowleft}

\def\id{{\mbf{\rm id}}}

\def\mbf#1{\hbox{\mathversion{bold}$#1$}} 
\def\modcomp{\:\tilde{\circ}\,} 

\def\openbull{\framebox[0.08in][c]{$\;$}} 

\def\re{{\mathbb R}} 

\def\shuffle{{\scriptscriptstyle \;\sqcup \hspace*{-0.05cm}\sqcup\;}}

\def\begals{\[\begin{aligned}}
\def\endals{\end{aligned}\]}
\def\begal{\begin{align*}}
\def\endal{\end{align*}}
\def\begce{\begin{center}}
\def\endce{\end{center}}
\def\begar{\begin{array}}
\def\endar{\end{array}}
\def\begeq{\begin{equation}}
\def\endeq{\end{equation}}
\def\begdi{\begin{displaymath}}
\def\enddi{\end{displaymath}}
\def\begdis{\begin{eqnarray*}}
\def\enddis{\end{eqnarray*}}
\def\begeqa{\begin{eqnarray}}
\def\endeqa{\end{eqnarray}}
\def\begdes{\begin{description}}
\def\enddes{\end{description}}
\def\begit{\begin{itemize}}
\def\endit{\end{itemize}}
\def\begen{\begin{enumerate}}
\def\enden{\end{enumerate}}
\def\beglar{\left[\begin{array}}
\def\endrar{\end{array}\right]}
\def\begle{\begin{lemma}}
\def\endle{\end{lemma}}
\def\begde{\begin{definition}}
\def\endde{\end{definition}}
\def\begth{\begin{theorem}}
\def\endth{\end{theorem}}
\def\begco{\begin{corollary}}
\def\endco{\end{corollary}}
\def\begprop{\begin{proposition}}
\def\endprop{\end{proposition}}	
\def\begex{\begin{example}}
\def\endex{\hfill\openbull \end{example}}
\def\begexer{\begin{exercise}}
\def\endexer{\end{exercise}}
\def\begalg{\begin{algo}}
\def\endalg{\end{algo}}
\def\begre{\noindent{\bf Remark}:\hspace*{0.05cm}}
\def\endre{\\}
\def\begres{\noindent{\bf Remarks}:\begin{enumerate}}
\def\endres{\end{enumerate} \par}
\def\begpr{\noindent{\em Proof:}$\;\;$}
\def\endpr{\hfill\bull}
\def\begtab{\begin{tabular}}
\def\endtab{\end{tabular}}
\def\rref#1{(\ref{#1})}
\newtheorem{lemma}{Lemma}[section]
\newtheorem{definition}{Definition}[section]
\newtheorem{theorem}{Theorem}[section]
\newtheorem{proposition}{Proposition}[section]
\newtheorem{corollary}{Corollary}[section]
\newtheorem{example}{Example}[section]
\newtheorem{algo}{Algorithm}[section]
\def\HAprod{\mbf{m}} 
\def\shuff#1#2{\mathbin{
		\hbox{\vbox{\hbox{\vrule \hskip#2 \vrule height#1 width 0pt}\hrule}\vbox{\hbox{\vrule \hskip#2 \vrule height#1 width 0pt\vrule }\hrule}}}}
\def\shuffl{{\mathchoice{\shuff{5pt}{3.5pt}}{\shuff{5pt}{3.5pt}}{\shuff{3pt}{2.6pt}}{\shuff{3pt}{2.6pt}}}}
\def\shuffle{{\, \shuffl \,}}

\def\xlongrightarrow#1{\longrightarrow}

\usepackage{forest}
\forestset{
  decor/.style = {
    label/.expanded = {[inner sep = 0.2ex, font=\unexpanded{\tiny}]right:{$#1$}}
  },
  root/.style = {minimum size = 0.1ex},
  decorated/.style = {
    for tree = {
      circle, fill, inner sep = 0.3ex, minimum size = 1.ex,
      grow' = south, l = 0, l sep = 1.2ex, s sep = 0.7em,
      fit = tight, parent anchor = center, child anchor = center,
      delay = {decor/.option = content, content =}
    }
  },
  default preamble = {decorated, root},
  begin draw/.code={\begin{tikzpicture}[baseline={([yshift=-0.5ex]current bounding box.center)}]},
}

\def\shuff#1#2{\mathbin{
      \hbox{\vbox{\hbox{\vrule \hskip#2 \vrule height#1 width 0pt}\hrule}\vbox{\hbox{\vrule \hskip#2 \vrule height#1 width 0pt\vrule }\hrule}}}}
\def\shuffl{{\mathchoice{\shuff{5pt}{3.5pt}}{\shuff{5pt}{3.5pt}}{\shuff{3pt}{2.6pt}}{\shuff{3pt}{2.6pt}}}}
\def\shuffle{{\, \shuffl \,}}

\begin{document}


\title[Post-Lie Structures in Affine Feedback]{On the Post-Lie Structure in\\[0.1cm] SISO Affine Feedback Control Systems}


\author[K.~Ebrahimi-Fard]{Kurusch Ebrahimi-Fard}
\address{Department of Mathematical Sciences, Norwegian University of Science and Technology (NTNU), 7491 Trondheim, Norway. Centre for Advanced Study (CAS), Drammensveien 78, 0271 Oslo, Norway.}
\email{kurusch.ebrahimi-fard@ntnu.no}
\urladdr{https://folk.ntnu.no/kurusche/}

\author[W.~S.~Gray]{W.~Steven~Gray}
\address{Department of Electrical and Computer Engineering, Old Dominion University, Norfolk, VA 23529 USA.}
\email{sgray@odu.edu}
\urladdr{http://www.ece.odu.edu/~gray/}

\author[Venkatesh~G.~S.]{Venkatesh~G.~S.}
\address{Department of Mathematical Sciences, Norwegian University of Science and Technology (NTNU), 7491 Trondheim, Norway.}
\email{subbarao.v.guggilam@ntnu.no}


\begin{abstract}
The main objective of this work is to show that the single-input, single-output (SISO) affine feedback group, a transformation group in the context of the affine feedback interconnection of Chen--Fliess series, is a post-group in the sense of Bai, Guo, Sheng and Tang.
\end{abstract}




\maketitle

\noindent
\tableofcontents

\thispagestyle{empty}


\section{Introduction}
\label{sect:intro}

\begin{figure}[tb]
\begin{center}	
		{\scalebox{0.8}{
		\begin{tikzpicture}[node distance= 20mm, bl/.style = {draw,rectangle},cross/.style={path picture={
				\draw[black]
				(path picture bounding box.south)
				-- (path picture bounding box.north) (path picture bounding box.west)
				-- (path picture bounding box.east);
		}}]
		\node[bl, thick, minimum width=3cm, minimum height=1.5cm] (Fc) {\Large${F_{\mathrm{f}}}$};
		\node[draw, thick, circle,minimum width=0.75 cm] (mult) [left of=Fc, xshift = -1cm] {\Large$+$};
		\node (v) [left of = mult, xshift = 0.5cm]{\Large$u$};		
		\node (vx) [below of = mult, yshift=0.9cm, xshift = -0.3cm]{\Large$v$};	
		\node[bl, thick, minimum width=3cm, minimum height=1.5cm] (Fd) [below of=Fc, yshift = -1cm] {\Large ${F_{\mathrm{fb}}}$};
		\node (y) [right of= Fc, xshift = 1.5 cm]{\Large$y$};
		\node (dum) [above of =mult]{};
		\draw[->, -latex, thick] (v) -- (mult) ;
		\draw[->, -latex, thick] (mult) -- node[midway, yshift = 0.25cm] (u){} (Fc);
		\draw[->, -latex, thick] (Fd) -| (mult) ;
		\draw[->, -latex, thick] (Fc) -- node[midway](conn){} (y)  ;
		\draw[->, -latex, thick] (conn.center) |- (Fd);
		\end{tikzpicture}}}
		\caption{Map $F_{\text{f}}$ in additive feedback with map $F_{\text{fb}}$.}
		\label{fig:additive-fb}
\end{center}
\end{figure}
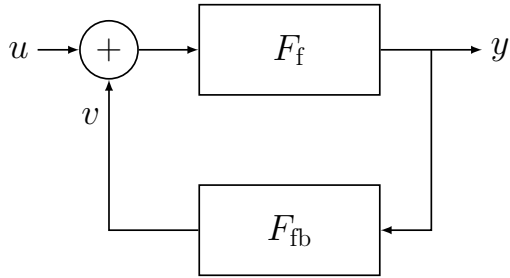

Consider two normed linear spaces of real-valued functions $U$ and $Y$ defined on a common interval $I$ and two maps $F_\text{f}:U \rightarrow Y$ and $F_{\text{fb}}:Y \rightarrow U$. An {\em additive feedback connection} with $F_\text{f}$ in the forward path and $F_{\text{fb}}$ in the feedback path is described by a
third map $F : U \rightarrow Y$ satisfying the feedback equations:
\begin{equation*}
\label{eq:FBeq}
y=F[u], \quad v=F_{\text{fb}}[y],\; \text{ and }\;  y=F_\text{f}[u+v]
\end{equation*}
for all $u \in U$. Such recursive interconnections are a central object of study in control theory \cite{Isidori_95}, where they are represented diagrammatically in Figure \ref{fig:additive-fb}.
The general objective is to determine $F_{\text{fb}}$ for a given $F_\text{f}$ so that $F$ has certain desirable properties such as having a finite operator norm.  An {\em affine feedback connection} is an important generalization of this set up where now the feedback equations take the form:
\begin{equation}
\label{eq:affFBeq}
y=F[u], \;
v_i=F_{\text{fb},i}[y],\ i=1,2,\  \text{and}\
y=F_\text{f}[uv_1+v_2]
\end{equation}
for all $u \in U$. Here the product $uv_1$ of real-valued functions is defined pointwise on $I$.		

It is well known that many types of feedback connections can be re-cast in terms of a feedback transformation group acting on a given map $F_\text{f}$ \cite{Brockett_78,Brockett_83,Brockett-Krishnaprasad_80,Gray-etal_SCL14,Gray-KEF_SIAM2017,Jakubczyk_90,Jakubczyk-Respondek_80,Tall-Respondek_03}.
The affine feedback transformation group presented in \cite{Gray-KEF_SIAM2017} describes the affine feedback connection when all the maps involved have Chen--Fliess series representations \cite{Fliess_81,Fliess_realizn_83,Isidori_95}. An important invariant under this group action is the {\em relative degree} of $F_\text{f}$. Loosely speaking, it describes the minimum number of integrations that $F_\text{f}$ applies to $u$. It can be defined and computed in terms of a smooth finite dimensional state space realization for $F_\text{f}$ when one is available \cite{Isidori_95,Nijmeijer-vanderSchaft_90} or in a
coordinate free manner when $F_\text{f}$ is a single-input, single-output (SISO) system and has only a Chen--Fliess series representation \cite{Gray-etal_AUTO14}.
Maps with relative degree are known to be left-invertible so that there exists an $F_{\text{fb}}$ rendering $F$ exactly linear \cite{Isidori_95,Nijmeijer-vanderSchaft_90}. So-called {\em feedback linearization} is a standard design method in control engineering for systems with nonlinear dynamics.
\medskip

Given $F_\text{f}$ and  $F_{\text{fb},i}$, solving the feedback equations \eqref{eq:affFBeq} for $F$ ultimately involves computing an inverse in the affine feedback transformation group. This calculation can be done inductively using the graded connected commutative polynomial Hopf algebra $\mathcal{H}$ of coordinate functions associated with the group \cite{Gray-KEF_SIAM2017}. Specifically, the generating series for the Chen--Fliess representation of $F$ is computed using the antipode of $\mathcal{H}$.
\medskip

Since $\mathcal{H}$ is a free commutative, right-handed Hopf algebra, there is an underlying right pre-Lie product defined on the space of infinitesimal characters associated with the group \cite{MenousPatras2015}. Its exact form was determined in  \cite{Gray-KEF_SIAM2017} by treating the group as a formal Lie group and computing the set of all left-invariant vector fields at the identity element. Subsequently, it was discovered by Foissy in \cite{Foissy_18} that the corresponding Lie algebra for this group coincides exactly with a suitably defined post-Lie algebra $\left(\mathfrak{g}_{\scriptscriptstyle{\mathcal{SISO}}}, \lhook, \left[\cdot,\cdot\right]\right)$ with $\left(\mathfrak{g}_{\scriptscriptstyle{\mathcal{SISO}}}, \llbracket\cdot,\cdot\rrbracket\right)$ being the derived Lie algebra. Since $\mathcal{H}$ is commutative, by the Cartier--Milnor--Moore theorem, the graded dual of $\mathcal{H}$, say $\mathcal{H}^{'}$, is isomorphic to the enveloping algebra $\mathcal{U}\left(\mathfrak{g}_{\mathcal{SISO}},\llbracket\cdot,\cdot\rrbracket\right)$ as connected graded Hopf algebras.

\medskip

It is shown in \cite{Bai2022postgroup} that the integration of a post-Lie algebra results in a post-group ~(see \cite{Alkaabi2023postgroup} for a detailed review of post-groups). The main contribution of the current work is to present a concise description of the post-group $\left(G,\cdot, \triangleleft\right)$ whose Grossman--Larson group $\left(G,\star\right)$ is the affine feedback transformation group.
The group $\left(G, \cdot \right)$, which is isomorphic to the semi-direct product of the additive group and the shuffle group of formal power series, went unnoticed in the previous works on affine feedback interconnections of Chen--Fliess series~\cite{Gray-KEF_SIAM2017,Foissy_18}. The post-group structure is derived purely from the \textit{ab initio} view of the affine feedback interconnection of Chen--Fliess series stemming from the earlier work of the first and second authors~\cite{Gray-KEF_SIAM2017}. Thus, the polynomial Hopf algebra of coordinate functions, $\mathcal{H}$, has a cointeracting bialgebra structure, where the two cointeracting coproducts $\delta$ and $\Delta$ come from dualizing the two group products $\left(G, \cdot\right)$ and the Grossman--Larson group product $\left(G, \star\right)$ associated with the post-group. The cointeraction picture in $\mathcal{H}$ is another contribution of the present work.

\medskip

{\bf{Acknowledgments}}: The first author is supported by the Research Council of Norway through project
302831 “Computational Dynamics and Stochastics on Manifolds” (CODYSMA). He also thanks the Centre for Advanced Study (CAS) in Oslo for support. The third author was supported by an ERCIM fellowship.

\medskip
\medskip

In the following, all algebras are defined on a field of characteristic zero (typically $\re$ if not specifically mentioned).


\section{Post-Groups: An Overview}
\label{sect:PG}

The notion of a post-group recently appeared in the work of Bai, Guo, Sheng and Tang \cite{Bai2022postgroup}. See \cite{Mencat2020,Mencat2021} for an equivalent formulation using the notion of a crossed morphism. Post-Lie algebras are the infinitesimal objects of post-groups as explained in Section \ref{sec:post-Lie-alg}. The main definition is given
below.

\begde \cite{Bai2022postgroup}
\label{de:Post--group}
\begin{enumerate}
\item[a)]
A (right) post-group is a group $(G,\cdot)$ equipped with a product $\triangleleft: G \times G \longrightarrow G$ such that:
\begen[(i)]
	\item For all $a \in G$, the right multiplication map $R^\triangleleft_{a} : G \longrightarrow G$
	$$
		R^\triangleleft_{a} (b) := b \triangleleft a
	$$
	is an automorphism of the group $\left(G, \cdot \right)$, that is,
	\begin{align*}
		R^\triangleleft_{a} (b \cdot c) = R^\triangleleft_{a} (b) \cdot R^\triangleleft_{a} (c).
	\end{align*}
	
\medskip

	\item For any $a,b \in G$, the following relation holds
	\begin{align}
	\label{eqn:weighted-associativity}
	R^\triangleleft_{a} \circ R^\triangleleft_{b} = R^\triangleleft_{R^\triangleleft_{a}(b)  \cdot a}.
	\end{align}
\enden

\item[b)] If $\left(G, \cdot \right)$ is Abelian, then $\left(G, \cdot, \triangleleft\right)$ is called a pre-group.
\end{enumerate}
\endde

\medskip

From \eqref{eqn:weighted-associativity}, one can deduce that the collection $\{R^\triangleleft_a : a \in G\}$ of automorphisms of the group $\left(G, \cdot \right)$ itself forms a group with the composition~(group) product defined by
\begin{align}
	\label{eqn:GLprod}
	b \star a := R^\triangleleft_{a} (b)  \cdot a = \left(b \triangleleft a\right) \cdot a
\end{align}
for all $a,b$ in the post-group $\left(G, \cdot, \triangleleft \right)$ such that $R^\triangleleft_a \circ R^\triangleleft_b = R^\triangleleft_{b \star a}$. Indeed, the following result confirms this statement (see \cite{Alkaabi2023postgroup} for a brief survey of the main results in \cite{Bai2022postgroup}).

\begin{theorem} \cite{Bai2022postgroup}
If $\left(G, \cdot, \triangleleft\right)$ is a post-group, then $\left(G, \star \right)$ forms a group called the {\em Grossman--Larson group} and shares the unit with $(G, \cdot)$. The inverse with respect to \eqref{eqn:GLprod} is
\begin{align*}
	a^{\star -1} =(R^\triangleleft_a)^{\circ -1} (a^{\cdot -1} ),
\end{align*}
and
\begin{align*}
\begin{aligned}
	R^\triangleleft : \left(G, \star \right)
	&\longrightarrow \mathrm{Aut}\left(\left(G, \cdot \right),\circ\right) \\
	a
	&\longmapsto R^\triangleleft_a
\end{aligned}
\end{align*}
is a group homomorphism.\footnote{Note that in \cite{Bai2022postgroup}, $\left(G, \star \right)$ was called the sub–adjacent group. The authors believe, however, that the name Grossman--Larson group \cite{Alkaabi2023postgroup} is more adequate in light of the origins of this notion in the work of Guin--Oudom \cite{GuinOudom2008preLie} on the Lie enveloping algebra of a pre-Lie algebra. See also \cite{MuntheKaas2015postLie} for the post-Lie case.}
\end{theorem}

\begre{
\label{rmk:postgroup}
\begin{enumerate}
\item Observe that for any $a,b \in G$
$$
	(a  \triangleleft b^{\star -1}) \star b = \big( (a  \triangleleft b^{\star -1})  \triangleleft b \big) \cdot b= a \cdot b.
$$

\item  Given a post-group $\left(G, \cdot, \triangleleft\right)$, the opposite post-group $\left(G, \centerdot, \blacktriangleleft\right)$  is defined in terms of the opposite group product, $a  \centerdot b := b \cdot a$, and for any $a, b \in G$ the product $a \blacktriangleleft b := a \cdot (a \triangleleft b) \cdot a^{\cdot -1}$ \cite{Alkaabi2023postgroup}.

\item
Natural examples of pre- and post-groups are found in the universal enveloping algebras of free pre-, respectively,  post-Lie algebras (see Section \ref{sec:post-Lie-alg}). In the former case, the pre-group is the Butcher group, i.e., the characters of the Butcher--Connes--Kreimer Hopf algebra of rooted trees \cite{GuinOudom2008preLie}. In the latter case, the post-group is the character group of the Munthe-Kaas--Wright Hopf algebra, known as Lie--Butcher group \cite{Alkaabi2023postgroup,MuntheKaas2015postLie}.
\end{enumerate}
}


\section{Preliminaries on Chen--Fliess Series}
\label{sect:FPS}

Consider the alphabet $X = \{x_0,x_1\}$. The free monoid over $X$ is denoted $\allwords$ and consists of the set of all non-commutative words over $X$, including the empty word which is denoted by $1$.  The subset of non-empty words is $X^{+} = X\setminus\{1\}$. The tensor algebra $T(X) := \bigoplus_{n \geq 0} \re X^{\otimes n}$ is isomorphic to the catenation algebra $\allpoly$ of non-commutative polynomials.

\begen[(i)]
\item $\left(\allpoly, \Delta_{\shuffle}, \epsilon\right)$ is the cocommutative bialgebra with the unshuffle coproduct defined by
\begin{align*}
	\Delta_{\shuffle}(x_i) = x_i \otimes 1 + 1 \otimes x_i
\end{align*}
and counit $\epsilon: \allpoly \to \re$
\begin{align*}
	\epsilon(w)  = 	\begin{cases}
					1 \quad \text{if} \; w = 1 \\
					0 \quad \text{otherwise}.
				\end{cases}
\end{align*}

\medskip

\begre{
Traditionally, the counit of $\allpoly$ is denoted by $\emptyset$ in the Chen--Fliess literature; the presentation also follows this custom, i.e., $\epsilon=\emptyset$.
}
\medskip

\item The convolution algebra of linear maps from the counital coalgebra $\left(\allpoly, \Delta_{\shuffle}, \emptyset\right)$ to $\re$, is given by the space of formal power series denoted by $\allseries{}$. The dual basis is given by $\{\emptyset\} \cup X^{+}$ such that $\eta\left(\xi\right) = 1$ if $\eta = \xi$ in $X^{+}$ and zero otherwise. An element $c \in \allseries{}$ is represented by
\begin{align*}
	c = c(1)\emptyset + \sum\limits_{\eta \in X^+} c(\eta) \eta.
\end{align*}
In the following, $\emptyset$ is not explicitly written unless needed.
\smallskip

The convolution product on $\allseries{}$ is the shuffle product, which is defined for all $c,d \in \allseries{}$ and $p \in \allpoly$ by
\begin{align*}
	\left(c \shuffle d\right)(p)
	&= m_{\re} \circ \left(c \otimes d\right) \circ \Delta_{\shuffle} (p).
\end{align*}
Here, $m_{\re}$ is the usual product in $\re$, and the unit element is the counit $\emptyset$.
\medskip

\item The maximal ideal of $\allseries{}$ are the so-called proper series, defined as
$$
	\allproperseries{} : = \{c \in \allseries{} : c(1) = 0\}.
$$
The group $G_{\shuffle} : = \left(1\emptyset + \allproperseries{}, \shuffle, \emptyset\right)$ is the normal subgroup of the unit group of $\allseries{}$, where
\begin{align*}
	1\emptyset + \allproperseries{} := \{c \in \allseries{}: c(1) = 1\} \subset  \allseries{}.
\end{align*}
\enden


\subsection{Chen--Fliess Series}
\label{sec:CF-series}

For $c \in \allseries{}$ and an integrable function $u$ (control signal), a Chen--Fliess series $F_c$ is a functional series such that the basis (words) of $\allseries{}$ are identified with iterated integrals of $u$. In the context of a dynamical system, a Chen--Fliess series expresses its input-output behavior and provides a coordinate-independent framework to analyze the intrinsic properties of the system.

\medskip

Let $u : [0,T] \longrightarrow \re$ be a $L_{1}$-measurable map over the compact interval $[0,T]$. Define the absolutely continuous function $U : [0,T] \longrightarrow \re$ such that $U(0) = 0$ as
\begin{align*}
	U(t) = \int_{0}^t d\tau u(\tau).
\end{align*}
Define a monoid morphism $F$ from the monoid of (dual) words $X^{+} \cup \{\emptyset\}$ to monoid of iterated integrals of $u$, where $F_\emptyset[u](t)=1$ is the constant function, and for all $\eta \in X^{+} \cup \{\emptyset\}$
\begin{align*}
	F_{x_1\eta}[u](t)
	&= \int_{0}^t dU_i(\tau)F_{\eta}[u](\tau),\\
	F_{x_0\eta}[u](t)
	&= \int_{0}^t d\tau F_{\eta}[u](\tau).
\end{align*}
Given a non-commutative formal power series $c \in \allseries{}$, the associated {\em Chen--Fliess series}, denoted by $F_{c}[u] : [0,T] \longrightarrow L_{\infty}([0,T]; \re)$, is a formal series defined as
\begin{align*}
	F_{c}[u](t) = c(1) + \sum_{\eta \in X^+} c(\eta)F_{\eta}[u](t).
\end{align*}
Observe that $F_{c \shuffle d}[u]=F_c[u]  F_d[u]$. For more details, see~\cite{Gray-Wang_SCL02,Thitsa-Gray_SIAM12,GS_thesis,Irina_thesis,Winter_Arboleda-etal_2015}.
\medskip

\begin{figure}[tb]
\includegraphics[scale= 0.8]{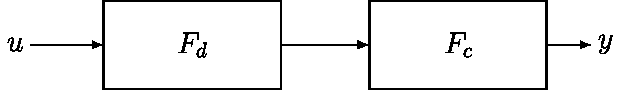}
\caption{Composition of Chen--Fliess series}
\label{fig:comp-CF-series}
\end{figure}

The composition of Chen--Fliess series
$$
	(F_c \circ F_d) [u]= F_c[F_d[u]],
$$
as illustrated in Figure~\ref{fig:comp-CF-series},
is defined in the obvious way, i.e., by making the function $F_d[u]$ the input of the Chen--Fliess series $F_c$. Using the multi-linearity of Chen's iterated integrals, this means looking at iterated integrals inside iterated integrals. Thanks to integration by parts the composition $F_c[F_d[u]]$ can be expressed as a bonafide Chen--Fliess series $F_{c \circ d}[u]$. The computationally challenging task is to present a transparent description of the formal power series $c \circ d \in \allseries{}$, called the {composition product} of $c$ and $d$ in $\in \allseries{}$, which is
\begin{align}
\label{eqn:comp-prod}
	c \circ d := c(1)\emptyset  + \sum_{\eta \in X^{+}} c\left(\eta\right) \eta \circ d.
\end{align}
It is computed in terms of its restriction to words, $\eta \circ d$, $\eta \in \allwords$, which is inductively defined in \cite{Ferfera_79,Ferfera_80,Gray-Li_05} by
\begin{align}
\label{eqn:comp-prod-ind}
	\emptyset \circ d
	&= \emptyset, \nonumber \\
	x_0\eta \circ d
	&= x_0\left(\eta \circ d\right), \\
	x_1\eta \circ d
	&= x_0\left(d \shuffle \left(\eta \circ d\right)\right) \nonumber.
\end{align}

\begre{~For an input-affine nonlinear state space model $\dot{z}(t) = f_0(z) + f_1(z)u$, $z(0) \in \re^n$ and $y= h(z)$,
where $f_{0},f_1 : \re^n \longrightarrow \re^n$ are analytic functions, and $h : \re^n \longrightarrow \re$ is also analytic, the {\em input-output} map $u \mapsto y$ is given {\em locally} by the Chen--Fliess series $y = F_c[u]$, where the coefficients are $c(1) = h(z(0))$ and $c(x_{i_1}x_{i_2}\cdots x_{i_m}) = L_{{i_m}}\cdots L_{{i_2}}L_{{i_1}}h (z(0))$ with $i_k \in \{0,1\}$ for all $k = 1,2,\ldots,m$. Here, $L_{{i_k}}$ denotes the Lie differential operator corresponding to vector field $f_{i_k}$ \cite{Fliess_realizn_83}.
}


\subsection{(Feed-forward) Chen--Fliess Series}

\begen
\item  Let $\allseries{}^2$ be the vector space of two-dimensional formal power series expressed in the basis $\{\text{e}_1,\text{e}_2\}$ using bold-face font
\begin{align*}
	\allseries{}^2 := \{\mathbf{c} = c_1\text{e}_1 + c_2\text{e}_2 : c_1,c_2 \in \allseries{}\}.
\end{align*}
Interchangeably, elements of $\allseries{}^2$ are represented using tuple notation
$$
	\mathbf{c} = [c_1 \;\; c_2]^t.
$$
\begre{
Note that the basis vectors $\text{e}_1$ and $\text{e}_2$ are referred to as $\delta$ and $1$ in \cite{Gray-KEF_SIAM2017}.
}

\medskip

\item For any $\mathbf{c} = [c_1 \;\; c_2]^t$ the corresponding (feed-forward) Chen--Fliess series is defined as
\begin{align*}
	F_{\mathbf{c}}[u] := uF_{c_1}[u] + F_{c_2}[u].
\end{align*}
The component $c_1$ is called the {\em multiplicative component} (or channel) and the component $c_2$ is termed the {\em additive component} (or channel).
\enden

\begin{figure}[tb]
\includegraphics[width = 0.65\textwidth, scale=0.8]{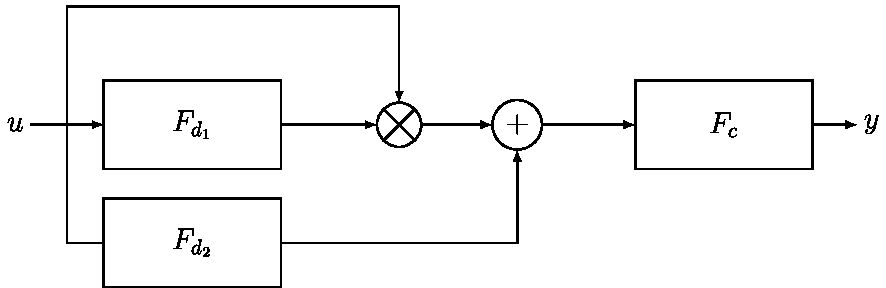}
\caption{Composition of Chen--Fliess series~$F_c[u]$ and $uF_{d_1}[u] + F_{d_2}[u]$.}
\label{fig:mixed_comp}
\end{figure}

The composition of Chen--Fliess series $F_{c} \circ F_{\mathbf{d}}$, where $c \in \allseries{}$ and $\mathbf{d} \in \allseries{}^2$ (illustrated in Figure~\ref{fig:mixed_comp}), is the Chen--Fliess series
$$
	F_{c \modcomp \mathbf{d}}:=F_{c} \circ F_{\mathbf{d}},
$$
corresponding to the formal power series $c \modcomp \mathbf{d}$ defined next.

\begde \cite[Def.~3.1]{Gray-KEF_SIAM2017}\label{def:mix-pro-single}
If $c \in \allseries{}$ and $\mathbf{d} = [d_1 \; d_2]^t \in \allseries{}^2$, then the mixed composition product
\begin{align}
\label{eqn:mixcomp-prod}
	c \modcomp \mathbf{d} := \sum_{\eta \in X^{*}} c(\eta) \eta \modcomp \mathbf{d} \in \allseries{}
\end{align}
is defined inductively on words $\eta \in \allwords$ by
\begin{align}
\label{eqn:mix-comp-ind}
	\emptyset \modcomp \mathbf{d}
	&= \emptyset,\nonumber\\
	x_0\eta \modcomp \mathbf{d}
	&= x_0(\eta \modcomp \mathbf{d}), \\
	x_1\eta \modcomp \mathbf{d}
	&= x_1(d_1 \shuffle (\eta \modcomp \mathbf{d})) + x_0(d_2 \shuffle (\eta \modcomp \mathbf{d})).\nonumber
\end{align}
\endde
\noindent Note that the mixed composition product \eqref{eqn:mix-comp-ind} is linear in its left argument (by definition).

\begth\label{th:dist-over-sh}\cite[Lem.~3.3]{Gray-KEF_SIAM2017}
If $c,d \in \allseries{}$ and $\mathbf{z} \in \allseries{}^2$, then
\begin{align}
\label{thm:distrib1}
	\left(c \shuffle d\right) \modcomp \mathbf{z}
	= \left(c \modcomp \mathbf{z}\right) \shuffle \left(d \modcomp \mathbf{z}\right).
\end{align}
\endth

\begth
\label{th:mix-assoc}
If $e,c \in \allseries{}$ and $\mathbf{d} \in \allseries{}^2$, then $e \circ \left(c \modcomp \mathbf{d}\right) = \left(e \circ c\right) \modcomp \mathbf{d}$.
\endth

\begpr
Considering the linearity of $\circ$ in its left argument in \rref{eqn:comp-prod}, it is sufficient to prove the following statement for words $\eta \in \allwords$, i.e., $\eta \circ \left(c \modcomp \mathbf{d}\right) = \left(\eta \circ c\right) \modcomp \mathbf{d}.$ The proof is by induction on length of the word $\eta$, here denoted by $\abs{\eta}$. The base case $\abs{\eta} = 0$ corresponds to $\eta = \emptyset$. By Definition~\ref{def:mix-pro-single} and \rref{eqn:comp-prod-ind}, $\emptyset \circ \left(c\modcomp \mathbf{d}\right) = \emptyset = \left(\emptyset \circ c\right) \modcomp \mathbf{d}.$ Assuming the hypothesis is true for all $\eta$ of length $\abs{\eta} = n$. There are two cases for a word of length $n+1$ to be considered, i.e.,~$x_0 \eta$ and $x_1 \eta$. Using \rref{eqn:comp-prod-ind},
\begin{align*}
	x_0\eta \circ \left(c\modcomp \mathbf{d}\right)
	&= x_0\left(\eta\circ\left(c \modcomp \mathbf{d}\right)\right) \\
	&= x_0\left(\left(\eta \circ c\right)\modcomp \mathbf{d}\right)\\
	&\stackrel{\rref{eqn:mix-comp-ind}}{=} \left(x_0\left(\eta \circ c\right)\right)\modcomp \mathbf{d}\\
	&= \left(x_0\eta \circ c\right) \modcomp \mathbf{d}.
\end{align*}
For the word $x_1\eta$, using \rref{eqn:comp-prod-ind} implies
\begin{align*}
	x_1\eta\circ\left(c \modcomp \mathbf{d}\right)
	&= x_0\left(\left(c \modcomp \mathbf{d}\right) \shuffle \left(\eta \circ \left(c \modcomp \mathbf{d} \right)\right)\right) \\
	&= x_0\left(\left(c \modcomp \mathbf{d}\right) \shuffle \left(\left(\eta \circ c\right)\modcomp \mathbf{d}\right)\right)\\
	&\stackrel{\eqref{thm:distrib1}}{=} x_0\left(\left(c  \shuffle \left(\eta \circ c\right)\right)\modcomp \mathbf{d}\right) \\
	&\stackrel{\rref{eqn:mix-comp-ind}}{=} \left(x_0\left(c \shuffle\left(\eta \circ c\right)\right)\right) \modcomp \mathbf{d}\\
	&\stackrel{\eqref{eqn:comp-prod-ind}}{=} \left(x_1\eta \circ c\right) \modcomp \mathbf{d}.
\end{align*}
Therefore, $e \circ \left(c \modcomp \mathbf{d}\right) = \left(e \circ c\right)\modcomp \mathbf{d}$.
\endpr

The mixed composition product \eqref{eqn:mixcomp-prod} is extended by applying the product component-wise. The extended product is denoted by $\triangleleft$ and defined next.

\begde \label{def:mix-pro}
If $\mathbf{c},\mathbf{d} \in \allseries{}^2$ with $\mathbf{c} = [c_1\;\; c_2]^t$, then
\begin{align}
	\triangleleft\,: \allseries{}^2 \times \allseries{}^2 &\longrightarrow \allseries{}^2 \nonumber \\
	(\mathbf{c},\mathbf{d}) &\longmapsto \left(\mathbf{c} \triangleleft \mathbf{d}\right)
	:= [c_1\modcomp \mathbf{d} \;\; c_2\modcomp \mathbf{d}]^t. \label{eqn:mix-pro}
\end{align}
\endde

\begre
It should be noted that given $\mathbf{c},\mathbf{d} \in \allseries{}^2$, the composition of (feed-forward) Chen--Fliess series
$$
	\left(F_{\mathbf{c}} \circ F_{\mathbf{d}}\right)[u] \neq F_{\mathbf{c} \triangleleft \mathbf{d}}[u].
$$
The $\triangleleft$ product is merely an extension of the mixed composition product \eqref{eqn:mixcomp-prod} component-wise in $\allseries{}^2$. The product encoding the composition $\left(F_{\mathbf{c}} \circ F_{\mathbf{d}}\right)[u]$ is presented in the next section.
\endre


\section{Affine Feedback Group: Grossman--Larson Group of a Post-Group}
\label{ssec:affine-feedb-grp}

Let $G := \{\mathbf{c} \in \allseries{}^2: c_1(1) = 1\}$ and thus, $G \cong \left(1 + \allproperseries{}\right) \times \allseries{}$ as sets.

\begth\label{th:odot-grp}
Let $\mathbf{c},\mathbf{d} \in G$ with $\mathbf{c} = [c_1 \;\; c_2]^t$ and $\mathbf{d} = [d_1 \;\; d_2]^t$. Define on $G$ the binary product
\begin{align}
\label{eqn:grp-prod-odot}
	\mathbf{c} \cdot \mathbf{d}
	:= 	\begin{bmatrix}
		c_1 \shuffle d_1 & c_2 + (c_1 \shuffle d_2)
		\end{bmatrix}^t.
\end{align}
Then $(G, \cdot)$ is a group with the identity element being
$$
	\mathbf{e}:=[1\emptyset \;\; 0]^t.
$$
The group inverse for $\mathbf{c} \in G$ is defined as
\begin{align*}
	\mathbf{c}^{\tiny{\cdot -1}} = 	
	\begin{bmatrix}
	c_1^{\shuffle -1} & -c_{1}^{\shuffle -1} \shuffle c_2
	\end{bmatrix}^t.
\end{align*}
\endth

\begpr
The definitions of the identity element and group inverse follow by straightforward verification. Associativity of \eqref{eqn:grp-prod-odot} is shown next. If $\mathbf{z} = [z_1 \;\; z_2]^t \in G$, then
\begin{align*}
	\left(\mathbf{c} \cdot \mathbf{d} \right) \cdot \mathbf{z}
	&= \left(\begin{bmatrix} c_1 \shuffle d_1 \\ c_2 + (c_1 \shuffle d_2) \end{bmatrix}\right) \cdot \mathbf{z} \\
	& = \begin{bmatrix} c_1\shuffle d_1 \shuffle z_1 \\
	c_2 + \left(c_1 \shuffle d_2 \right) + \left(c_1 \shuffle d_1 \shuffle z_2\right) \end{bmatrix}
	= \begin{bmatrix} c_1 \shuffle \left(d_1 \shuffle z_1\right) \\
	c_2 + c_1 \shuffle \left(d_2 + d_1 \shuffle z_2 \right)\end{bmatrix}
	= \mathbf{c} \cdot (\mathbf{d} \cdot \mathbf{z}).
\end{align*}
\endpr

The group $(G, \cdot)$ is isomorphic to the semi-direct product of $G_{\shuffle}$ and the additive group of formal power series,
\begin{align*}
	(G, \cdot) \cong G_{\shuffle} \ltimes (\allseries{},+),
\end{align*}
with the split exact sequence of groups given by
\begin{displaymath}
\begin{tikzcd}
	0 \arrow{r} & \allseries{} \arrow{r}{i} & G \arrow{r}{\pi_1} & G_{\shuffle} \arrow{r} & 1 \, ,
\end{tikzcd}
\end{displaymath}
where $i$ is the monomorphism $i(c) = 1\text{e}_1 + c\text{e}_2$ for all $ c \in \allseries{}$, and
$\pi_1$ is the canonical projection onto the first component.

\smallskip

\begre{~$i(\allseries{}) = \{1\text{e}_1 + c\text{e}_2 : c \in \allseries{}\}$ and $G_{\shuffle}$ are the underlying groups in the definition of the pre-groups corresponding to the (SISO) additive feedback and multiplicative feedback interconnections of Chen--Fliess series, respectively.}

\begth \label{thm:dist-over-odot}
If $\mathbf{x}, \mathbf{y}, \mathbf{d}\in G$ and $\mathbf{e}=[1\emptyset \;\; 0]^t$ is the identity element of $\left(G, \cdot \right)$, then
\begin{enumerate}[(i)]
\item $\left(\mathbf{x} \cdot \mathbf{y}\right) \triangleleft \mathbf{d} = \left(\mathbf{x} \triangleleft \mathbf{d}\right) \cdot \left(\mathbf{y} \triangleleft \mathbf{d}\right)$
\medskip

\item $\mathbf{d} \triangleleft \mathbf{e} = \mathbf{d}$
\medskip

\item $\mathbf{e} \triangleleft \mathbf{d} = \mathbf{e}$.
\end{enumerate}
\endth

\begpr Statements (ii) and (iii) follow from a straightforward calculation using Definition~\ref{def:mix-pro}.

(i) If $\mathbf{x} = x_1 \text{e}_1 + x_2 \text{e}_2$ and $\mathbf{y} = y_1 \text{e}_1 + y_2 \text{e}_2$, then using~\rref{eqn:grp-prod-odot} and Definition~\ref{def:mix-pro},
\begin{align*}
	\left(\mathbf{x} \cdot\mathbf{y}\right) \triangleleft  \mathbf{d}
	&= \begin{bmatrix} x_1 \shuffle y_1 \\ x_2 + \left(y_2 \shuffle x_1\right)
	\end{bmatrix} \triangleleft  \mathbf{d} \\
	&= \begin{bmatrix}
	\left(x_1 \shuffle y_1\right) \modcomp \mathbf{d} \\
	\left(x_2 + \left(x_1 \shuffle y_2 \right)\right) \modcomp \mathbf{d}
	\end{bmatrix}\\
	&\stackrel{\eqref{thm:distrib1}}{=} \begin{bmatrix}
	\left(x_1 \modcomp \mathbf{d}\right) \shuffle \left(y_1 \modcomp \mathbf{d}\right) \\
	\left(x_2 \modcomp \mathbf{d}\right) + \left( \left(x_1 \modcomp \mathbf{d}\right)\shuffle \left(y_2\modcomp \mathbf{d}\right)\right)
	\end{bmatrix} \\
	&\stackrel{\rref{eqn:grp-prod-odot}}{=} \begin{bmatrix}
	x_1 \modcomp \mathbf{d} \\ x_2 \modcomp \mathbf{d}
	\end{bmatrix} \cdot \begin{bmatrix}
	y_1 \modcomp \mathbf{d} \\ y_2 \modcomp \mathbf{d}
	\end{bmatrix}
	= \left(\mathbf{x} \triangleleft \mathbf{d}\right) \cdot \left(\mathbf{y} \triangleleft \mathbf{d}\right).
\end{align*}
\endpr

\begth
If $\mathbf{c}, \mathbf{d}, \mathbf{h}\in G$, then
\begin{align*}
	\left(\mathbf{c} \triangleleft \mathbf{d}\right) \triangleleft \mathbf{h}
	&= \mathbf{c} \triangleleft \left(\left(\mathbf{d} \triangleleft \mathbf{h}\right) \cdot \mathbf{h}\right).
\end{align*}
\endth

\begpr In light of the component-wise structure involved, it is sufficient to show that for all $c \in \allseries{}$
\begin{align*}
	\left(c \modcomp \mathbf{d}\right) \modcomp \mathbf{h}
	= c \modcomp \left(\left(\mathbf{d}\triangleleft \mathbf{h}\right)\cdot \mathbf{h}\right).
\end{align*}
In addition, since $\modcomp$ is linear in its left argument, one only needs to prove that for all words $\eta \in \allwords$:
\begin{align*}
	\left(\eta \modcomp \mathbf{d}\right) \modcomp \mathbf{h}
	= \eta \modcomp \left(\left(\mathbf{d}\triangleleft \mathbf{h}\right)\cdot \mathbf{h}\right).
\end{align*}
This is shown by induction on the length of $\eta$. The base case $\eta = \emptyset$ is trivial. Assume that the hypothesis is true for all $\eta = x_{i_1}x_{i_2}\cdots x_{i_k}$, where $i_j \in\{ 0,1\}$ for $j = 1,2,\ldots,k$. There are two cases to consider:

\begen[(i)]
\item
\begin{align*}
	\left(x_0\eta \modcomp \mathbf{d}\right) \modcomp \mathbf{h}
	&= \left(x_0\left(\eta \modcomp \mathbf{d}\right)\right) \modcomp \mathbf{h}\\
	&= x_0\left(\left(\eta \modcomp \mathbf{d}\right) \modcomp \mathbf{h}\right)\\
	&= x_0\left(\eta \modcomp \left(\left(\mathbf{d} \triangleleft \mathbf{h}\right) \cdot \mathbf{h}\right)\right) \\
	&= x_0\eta \modcomp \left(\left(\mathbf{d} \triangleleft \mathbf{h}\right) \cdot \mathbf{h}\right)
\end{align*}

\item
\begin{align*}
	\left(x_1\eta \modcomp \mathbf{d}\right) \modcomp \mathbf{h}
	&= \left[ x_1\left(d_1 \shuffle \left(\eta \modcomp \mathbf{d}\right)\right)
	+ x_0\left(d_2 \shuffle \left(\eta \modcomp \mathbf{d}\right)\right)\right] \modcomp \mathbf{h} \\
	&= x_1 \left(h_1 \shuffle \left(d_1 \modcomp \mathbf{h}\right) \shuffle \left(\left(\eta \modcomp \mathbf{d}\right)\modcomp \mathbf{h}\right) \right)\\
	&\quad \quad + x_0 \left(\left(\left(h_2 \shuffle \left(d_1 \modcomp \mathbf{h}\right)\right)
		+ \left(d_2 \modcomp \mathbf{h}\right)\right) \shuffle \left(\left(\eta \modcomp \mathbf{d}\right)\modcomp \mathbf{h} \right)\right) \\
	&= x_1\left(\left(\left(\mathbf{d}\triangleleft \mathbf{h}\right)\cdot \mathbf{h}\right)_1
			\shuffle \left(\eta \modcomp \left(\left(\mathbf{d} \triangleleft \mathbf{h}\right) \cdot \mathbf{h}\right)\right)\right) \\
	&\quad \quad + x_0\left(\left(\left(\mathbf{d}\triangleleft \mathbf{h}\right)\cdot \mathbf{h}\right)_2
			\shuffle \left(\eta \modcomp \left(\left(\mathbf{d} \triangleleft \mathbf{h}\right)\cdot \mathbf{h}\right)\right)\right)\\
	&= x_1\eta \modcomp \left(\left(\mathbf{d} \triangleleft \mathbf{h}\right)\cdot\mathbf{h}\right).
\end{align*}
\enden
\endpr

Lemma~$3.6$~(statement ($2$)) and Theorem~$3.8$ in \cite{Gray-KEF_SIAM2017} assert that the map
\begin{align*}
	R^{\triangleleft}_{\mathbf{c}} : G &\longrightarrow G \\
	\mathbf{b} &\longmapsto \mathbf{b} \triangleleft \mathbf{c}
\end{align*}
is a bijection on $G$. Thus, Theorem~\ref{th:odot-grp} and Theorem~\ref{thm:dist-over-odot} imply directly the following result.

\begth\label{th:affine_feedb}
$\left(G, \cdot, \triangleleft \right)$ is a post-group with Grossman--Larson group $\left(G, \star \right)$ and Grossman--Larson product defined by
\begin{align*}
	\star : G \times G &\longrightarrow G \\
	\left(\mathbf{c},\mathbf{d}\right) &\longmapsto \left(\mathbf{c} \triangleleft \mathbf{d}\right) \cdot \mathbf{d}.
\end{align*}
In particular, if $\mathbf{c} = c_1 \mathrm{e}_1 + c_2 \mathrm{e}_2$ and $\mathbf{d} = d_1 \mathrm{e}_1 + d_2 \mathrm{e}_2$, then
\begeq
\label{eqn:grp-prod-affine-cord}
	\mathbf{c} \star \mathbf{d}
	= \begin{bmatrix}
		\left(c_1 \modcomp \mathbf{d}\right) \shuffle d_1 & c_2\modcomp \mathbf{d}
						+ \left(c_{1}\modcomp \mathbf{d}\right) \shuffle d_2
	   \end{bmatrix}^t.
\endeq
\endth

\begin{figure}[tb]
\includegraphics[width = 0.45\textwidth, scale=0.6]{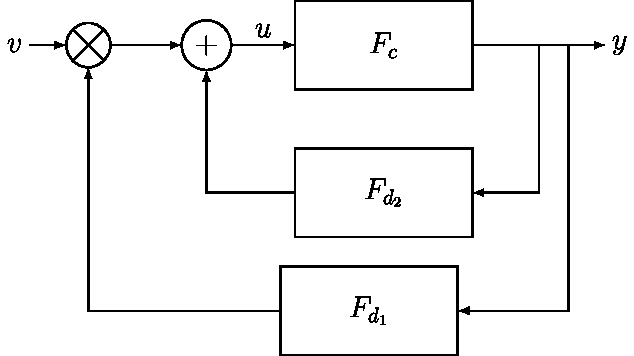}
\caption{Affine Feedback Interconnection of $F_{c}$ with $F_{\mathbf{d}}$}
\label{fig:affine-feedback2}
\end{figure}

\noindent Observe that~\rref{eqn:grp-prod-affine-cord} and Definition~$3.5$ in~\cite{Gray-KEF_SIAM2017} are identical. Thus, the Grossman--Larson group of the post-group $\left(G,\cdot, \triangleleft \right)$ is the {\em affine feedback group}. It is the transformation group related to the affine feedback interconnection of Chen--Fliess series depicted in Figure~\ref{fig:affine-feedback2}.

\section{From post-Lie Group to post-Lie Algebra}
\label{sec:post-Lie-alg}

Let $\mathbb{K}$ denote the base field of characteristic zero. Given a binary product $\diamond$ on a $\mathbb{K}$-vector space $V$, the associator is defined for elements $x,y,z \in V$ as
$$
	{\mathrm{a}}_{\diamond}(x,y,z)
	:= (x \diamond y) \diamond z - x \diamond (y \diamond z).
$$
When different from zero, the associator expresses the non-associativity of the product $\diamond$.

\begin{definition} \cite{MKW2008, BV2007}
\label{def:postLie}
A (right) \emph{post-Lie algebra} $(\mathfrak{g}, [ -,-], \lhook)$ consists of a Lie algebra $\mathfrak g$ with Lie bracket $[ -,-]$ and a product $\lhook : {\mathfrak g} \otimes {\mathfrak g} \rightarrow \mathfrak g$, called the post-Lie product, such that the following relations hold for all elements $x,y,z \in \mathfrak g$
\begin{align*}
	[x,y] \lhook z &= [x \lhook z , y] + [x , y \lhook z],\\
	z \lhook [x,y] &= {\mathrm{a}}_{\lhook }(z,x,y) - {\mathrm{a}}_{\lhook }(z,y,x).
\end{align*}
\end{definition}

Given a post-Lie algebra, a second Lie bracket can be defined on ${\mathfrak g}$
\begin{equation}
\label{postLie3}
	\llbracket x,y \rrbracket := x \lhook y - y \lhook x + [x,y], \quad \forall x,y \in {\mathfrak g}.
\end{equation}

If the Lie algebra $\mathfrak g$ in Definition \ref{def:postLie} has a trivial Lie bracket, then the post-Lie algebra reduces to a (right) {\textit{pre-Lie algebra}} $(\mathfrak{g}, \lhook)$, defined in terms of the (right) pre-Lie identity
\begin{equation}
\label{preLie1}
	(x \lhook y) \lhook z - x \lhook (y \lhook z)
	= (x \lhook z) \lhook y - x \lhook (z \lhook y).
\end{equation}
In this case, the Lie bracket in \eqref{postLie3} reduces to
\begin{equation}
\label{preLie2}
	\llbracket x,y \rrbracket := x \lhook y - y \lhook x.
\end{equation}
The pre-Lie relation \eqref{preLie1} implies that \eqref{preLie2} is a Lie algebra. Hence, anti-symmetrization is a covariant functor from the category of pre-Lie algebras over $\mathbb{K}$ to Lie algebras over $\mathbb{K}$. For more details on pre--Lie algebras, refer\cite{Manchon_11}.

\begin{example}
Consider the associative $\re$-Weyl algebra $\re[x_1,\ldots,x_n,\partial_1,\ldots,\partial_n]$ with the relations $\partial^j. x_j - x_j. \partial^j = 1$ and $\partial^j\partial^i-\partial^i\partial^j=0$. Let $\mathcal{W}$ be the $\re$-span of
$W = \{x_1^{i_1}x_2^{i_2}\cdots x_n^{i_n}\partial_i: i_1,i_2,\ldots, i_n \geq 0\}$.
Following \cite{Burde2006}, $\left(\mathcal{W},\blacktriangleleft\right)$ is a right pre-Lie algebra with pre-Lie product
\begin{align*}
\left(x_1^{i_1}x_2^{i_2}\cdots x_n^{i_n}\partial_i \right) \blacktriangleleft \left(x_1^{j_1}x_2^{j_2}\cdots x_n^{j_n}\partial_j\right) &=  i_j x_1^{i_1+j_1}x_2^{i_2+j_2}\cdots x_j^{i_j+j_i-1}\cdots x_n^{i_n+j_n} \partial_i.
\end{align*}
Dropping the commutativity among the differential operators, $\partial_1, \partial_2 \ldots, \partial_n$ results in a post-Lie algebra structure on $\mathcal{W}$. The scenario corresponds to polynomial vector fields on manifolds with constant torsion and zero curvature (see \cite{MuntheKaas2015postLie} for more details).
\end{example}

Given a post-Lie group $\left(G,\cdot, \triangleleft\right)$, let $\left(\mathfrak{g},\left[-,-\right]\right)$ be the Lie algebra of $\left(G, \cdot\right)$. Let the differential of the product $\triangleleft$ (at the identity) be denoted by $\lhook$. Proposition~$4.7$ of~\cite{Bai2022postgroup} states that $\left(\mathfrak{g},\left[-,-\right], \lhook\right)$ forms a post-Lie algebra. Furthermore, the derived Lie algebra $\left(\mathfrak{g}, \llbracket-,-\rrbracket\right)$ is the Lie algebra of the Grossman--Larson group of the post-Lie group.


\subsection{Post-Lie Algebra of the Affine Feedback}
\label{ssec:affine-post-Lie}

The Lie algebra (tangent space) $T_{\mathbf{e}}G$ of the (formal) Lie group $\left(G,\cdot\right)$ is denoted by $\left(\mathfrak{g}, \left[-,-\right]\right)$. As a vector space, $\mathfrak{g}$ is isomorphic to the Cartesian product of formal power series $\allproperseries{} \times \allseries{}$ and is spanned by the basis $\emptyset \text{e}_2$, $x_{i_1}x_{i_2}\cdots x_{i_k}\text{e}_j$ for all $k \geq 1$, $j \in \{1,2\}$ and
$i_{\ell} \in \{ 0,1\}$, $\ell = 1,2,\ldots,k$.\footnote{In \cite{Foissy_18}, $\mathfrak{g}$ was restricted to the space of polynomials viz. $\allpoly{} \times \allpoly{}$.}
Equivalently, $\mathfrak{g} \cong \allseries{}^2/(\emptyset e_1)$ as vector spaces.

\medskip

Let $\eta,\zeta \in X^{+} \cup \{\emptyset\}$ and $s,t \in \re$. Then $\exp\left(t\eta\text{e}_i\right)$ and $\exp\left(s\zeta\text{e}_j\right)$ are the group elements in the post--group $G$. The exponentials are defined by the respective associative products defined in the universal enveloping algebra of $\mathfrak{g}$~\cite{MuntheKaas2015postLie}. The product remains inconsequential for understanding the derivation of the post--Lie algebra associated with the post--group $\left(G,\cdot,\triangleleft\right)$, as only terms upto first order are used. Thus,
\begin{align*}
\exp(t\eta\text{e}_i) &= \mathbf{e} + t\eta\text{e}_i + \mathfrak{o}(t) \\
\exp(s\zeta\text{e}_j) &= \mathbf{e} + s\zeta\text{e}_j + \mathfrak{o}(s)
\end{align*}
where $\mathfrak{o}(t) \in t^2\re_p\langle X\rangle[[t]]$.
\medskip

The Lie bracket underlying the group $\left(G,\cdot\right)$ is computed next, which is the Lie bracket in the defintion of post--Lie algebra associated with $(G,\cdot, \triangleleft)$. For more details on deriving the post--Lie algebra associated with the post--group, refer~\cite{Bai2022postgroup}. 
\medskip

Let $\eta,\zeta \in X^{+} \cup \{\emptyset\}$ and $s,t \in \re$.
\begin{align*}
\exp(t\eta\text{e}_{i})\cdot \exp(s\zeta\text{e}_j) &= \left(\mathbf{e} + t\eta \text{e}_i + \mathfrak{o}(t)\right).\left(\mathbf{e} + s \zeta \text{e}_j + \mathfrak{o}(s)\right)\\
&\stackrel{\eqref{eqn:grp-prod-odot}}{=}\mathbf{e}
				+ t\eta \text{e}_i
				+ s \zeta \text{e}_j
				+ st\left(\eta\shuffle \zeta\right)\text{e}_j\delta_{i,1} + \mathfrak{o}(st),
\end{align*}
where the Kronecker delta is defined by $\delta_{i,1} := 1$ if $i = 1$ and zero otherwise. Hence, the bracket $[-,-]$ of the Lie algebra $\mathfrak{g}$ is
\begin{align}\label{eq:Lie-bracket-def}
	\left[ \eta  \text{e}_i, \zeta  \text{e}_j\right]
	&= \left(\eta\shuffle \zeta\right)\text{e}_j \delta_{i,1} - \left(\eta \shuffle \zeta\right)\text{e}_i \delta_{j,1}, \nonumber\\
	&=\begin{cases}
		\left(\eta \shuffle \zeta\right)\text{e}_2 	&\quad \text{if}\; i = 1, j = 2 \\
		-\left(\eta \shuffle \zeta\right)\text{e}_2 	&\quad \text{if}\; i = 2, j = 1. \\
		0 								&\quad \text{if}\; i = j.
	\end{cases}
\end{align}
\medskip

The differential of the action map $\triangleleft : G \times G \longrightarrow G$ is the map $\lhook : \mathfrak{g} \otimes \mathfrak{g} \longrightarrow \mathfrak{g}$, where for all $t,s \in \re$ and $\mathrm{u},\mathrm{v} \in \mathfrak{g}$:
\begin{align}
	\label{eqn:mixcomp-lin}
\exp(t\mathrm{u}) \triangleleft \exp(s\mathrm{v}) &= \left(\mathbf{e}
	+ t\mathrm{u} + \mathfrak{o}(t)\right) \triangleleft \left(\mathbf{e}
	+ s\mathrm{v} + \mathfrak{o}(s)\right) \nonumber \\
	&= \left(\mathbf{e} \triangleleft \left(\mathbf{e}
	+ s\mathrm{v} + \mathfrak{o}(s)\right)\right)  
	+ \left(t\mathrm{u} \triangleleft \left(\mathbf{e}
	+ s\mathrm{v} + \mathfrak{o}(s)\right)\right) + \left(\mathfrak{o}(t) \triangleleft \left(\mathbf{e}
	+ s\mathrm{v} + \mathfrak{o}(s)\right) \right) \nonumber \\
&=	\mathbf{e} + t\mathrm{u} + st(\mathrm{u} \lhook \mathrm{v}) + \mathfrak{o}(t) + \mathfrak{o}(st).
\end{align}
\medskip

The computation of the product $\lhook$ for any two Lie algebra elements is defined inductively next. 

\begth\label{th:mixpro-lin-ind}\cite{Foissy_18}
If $\eta, \zeta \in X^{+} \cup \{\emptyset\}$, then for $i,j = 1,2$
\begin{subequations}
\begin{align}
	\emptyset \mathrm{e}_i \lhook \zeta \mathrm{e}_j
	&= 0											\label{subeq1:emp} \\
	x_0\eta \mathrm{e}_i \lhook \zeta \mathrm{e}_j
	&= x_0\left(\eta \mathrm{e}_i \lhook \zeta \mathrm{e}_j\right) 	\label{subeq2:x_0} \\
	x_1\eta \mathrm{e}_i \lhook \zeta \mathrm{e}_1
	&= x_1\left(\eta \mathrm{e}_i \lhook \zeta \mathrm{e}_1\right)
		+ x_1\left(\eta \shuffle \zeta\right)\mathrm{e}_i			\label{subeq3:x_1e_1} \\
	x_1\eta \mathrm{e}_i \lhook \zeta \mathrm{e}_2
	&= x_1\left(\eta \mathrm{e}_i \lhook \zeta \mathrm{e}_2\right)
		+ x_0\left(\eta \shuffle \zeta\right)\mathrm{e}_i. 			\label{subeq3:x_1e_2}
\end{align}
\end{subequations}
\endth

\begpr Let $t,s \in \re$.
\begen[(i)]
\item Equality \eqref{subeq1:emp} follows from
\begin{align*}
	\exp\left(t\emptyset\mathrm{e}_i\right)\triangleleft \exp\left(s\zeta\mathrm{e}_j\right)&= \left(\mathbf{e}
	+ t\emptyset \mathrm{e}_2 + \mathfrak{o}(t)\right) \triangleleft \left(\mathbf{e}
	+ s\zeta \mathrm{e}_j + \mathfrak{o}(s)\right) \\
	&\stackrel{\eqref{eqn:mix-comp-ind}}{=} \mathbf{e} + t \emptyset \mathrm{e}_2.
\end{align*}

\item Equality \eqref{subeq2:x_0} follows from
\begin{align*}
\exp\left(tx_0\eta\mathrm{e}_i\right) \triangleleft \exp\left(s\zeta\mathrm{e}_j\right) &= \left(\mathbf{e}
		+ t x_0\eta \mathrm{e}_i + \mathfrak{o}(t)\right) \triangleleft \left(\mathbf{e} + s\zeta \mathrm{e}_j + \mathfrak{o}(s)\right) \\
		&= \mathbf{e} + t \left[x_0\eta \mathrm{e}_i \triangleleft \left(\mathbf{e} + s\zeta \mathrm{e}_j + \mathfrak{o}(s)\right)\right] + \mathfrak{o}(t) + \mathfrak{o}(st)\\
		&\stackrel{\eqref{eqn:mix-comp-ind}}{=} \mathbf{e}
			+ t x_0\left[\eta \mathrm{e}_i \triangleleft \left(\mathbf{e} + s\zeta \mathrm{e}_j + \mathfrak{o}(s)\right)\right] + \mathfrak{o}(t) + \mathfrak{o}(st) \\
		&\stackrel{\eqref{eqn:mixcomp-lin}}{=} \mathbf{e}
			+ t x_0\left[\eta \mathrm{e}_i + s\left(\eta \mathrm{e}_i \lhook \zeta \mathrm{e}_j\right) + \mathfrak{o}(s)\right] + \mathfrak{o}(t) + \mathfrak{o}(st)\\
		&= \mathbf{e} + tx_0\eta \mathrm{e}_i + st x_0\left(\eta \mathrm{e}_i \lhook \zeta \mathrm{e}_j\right) + \mathfrak{o}(t) + \mathfrak{o}(st).
\end{align*}

\item Equality \eqref{subeq3:x_1e_1} follows from
\begin{align*}
\exp\left(tx_1\eta\mathrm{e}_i\right) \triangleleft \exp\left(s\zeta\mathrm{e}_1\right) &=
	\left(\mathbf{e}+ tx_1 \eta \mathrm{e}_i + \mathfrak{o}(t)\right) \triangleleft \left(\mathbf{e} + s \zeta \mathrm{e}_1 + \mathfrak{o}(s)\right) \\
	&= \mathbf{e} + t\left[x_1\eta \mathrm{e}_i \triangleleft \left(\mathbf{e} + s \zeta \mathrm{e}_1 + \mathfrak{o}(s)\right)\right] + \mathfrak{o}(t) +\mathfrak{o}(st) \\
	&\stackrel{\eqref{eqn:mix-comp-ind}}{=} \mathbf{e}
	+ t x_1\left[ \left(1 + s\zeta\right) \shuffle \left(\eta \modcomp \left(\mathbf{e}
	+ s \zeta \mathrm{e}_1 + \mathfrak{o}(s)\right)\right)\right]\mathrm{e}_i + \mathfrak{o}(t) + \mathfrak{o}(st)\\
	&\stackrel{\eqref{eqn:mix-pro}}{=} \mathbf{e} + t x_1 [\left(\eta \mathrm{e}_i \triangleleft \left(\mathbf{e}
	+ s\zeta \mathrm{e}_1 + \mathfrak{o}(s)\right)\right) + \\
&\quad \quad \left(s \zeta \shuffle \left(\eta \modcomp \left(\mathbf{e}
	+ s\zeta \mathrm{e}_1 + \mathfrak{o}(s)\right)\right)\right)\mathrm{e}_i] 
+ \mathfrak{o}(t) + \mathfrak{o}(st)\\
	&\stackrel{\eqref{eqn:mixcomp-lin}}{=} \mathbf{e} + tx_1\left[\eta \mathrm{e}_i
	+ s \left(\eta \mathrm{e}_i \lhook \zeta \mathrm{e}_1\right) + \left(s\zeta \shuffle \eta\right)\mathrm{e}_i + \mathfrak{o}(s)\right] + \mathfrak{o}(t) + \mathfrak{o}(st) \\
	&= \mathbf{e} + tx_1\eta \mathrm{e}_i
	+ st \left[x_1\left(\eta \mathrm{e}_i \lhook \zeta \mathrm{e}_1\right)	+ x_1\left(\eta \shuffle \zeta \right)\mathrm{e}_i\right] + \mathfrak{o}(t) + \mathfrak{o}(st).
\end{align*}
\medskip

\item Equality \eqref{subeq3:x_1e_2} follows from
\begin{align*}
\exp\left(tx_1\eta\mathrm{e}_i\right) \triangleleft \exp\left(s\zeta\mathrm{e}_1\right) &=
	\left(\mathbf{e} + tx_1\eta \mathrm{e}_i + \mathfrak{o}(t)\right) \triangleleft \left(\mathbf{e} + s \zeta \mathrm{e}_2 + \mathfrak{o}(s)\right)\\
	&= \mathbf{e} + t\left[x_1\eta \mathrm{e}_i \triangleleft \left(\mathbf{e} + s \zeta \mathrm{e}_2 + \mathfrak{o}(s)\right)\right] + \mathfrak{o}(t) + \mathfrak{o}(st) \\
	&\stackrel{\eqref{eqn:mix-comp-ind}}{=} \mathbf{e} + tx_1\left[\eta \modcomp \left(\mathbf{e}
		+ s\zeta \mathrm{e}_2 + \mathfrak{o}(s)\right)\right]\mathrm{e}_i
		 \\
&\quad \quad + tx_0\left[s\zeta \shuffle \left(\eta \modcomp \left(\mathbf{e}
		+ s\zeta \mathrm{e}_2 + \mathfrak{o}(s)\right)\right)\right]\mathrm{e}_i + \mathfrak{o}(t) + \mathfrak{o}(st) \\
	&\stackrel{\eqref{eqn:mix-pro}}{=} \mathbf{e} + tx_1\left(\eta \mathrm{e}_i \triangleleft \left(\mathbf{e}
		+ s\zeta \mathrm{e}_2 + \mathfrak{o}(s)\right)\right)
		+ tx_0\left[\left(s\zeta \shuffle \eta\right)\mathrm{e}_i + \mathfrak{o}(s)\right]  \\
& \quad \quad \quad + \mathfrak{o}(t) + \mathfrak{o}(st) \\
&\stackrel{\eqref{eqn:mixcomp-lin}}{=} \mathbf{e} + tx_1\left(\eta \mathrm{e}_i
	+ s\left(\eta \mathrm{e}_i\lhook \zeta \mathrm{e}_i\right) + \mathfrak{o}(s)\right) + tx_0\left[\left(s\zeta \shuffle \eta\right)\mathrm{e}_i + \mathfrak{o}(s)\right] \\
& \quad \quad \quad + \mathfrak{o}(t) + \mathfrak{o}(st) \\
	&= \mathbf{e} + tx_1\eta \mathrm{e}_i
	+ st \left[x_1\left(\eta \mathrm{e}_i \lhook \zeta \mathrm{e}_2\right)
	+ x_0\left(\eta \shuffle \zeta\right)\mathrm{e}_i\right] + \mathfrak{o}(t) + \mathfrak{o}(st).
\end{align*}
\enden
\endpr
\medskip

Following Proposition~$4.7$ of \cite{Bai2022postgroup} gives the next theorem, which was first observed in \cite{Foissy_18} by an entirely different approach.

\begth\label{th:derived-Lie-affine-feedb}\cite{Foissy_18}
The triple
$\left(\mathfrak{g}, \left[\cdot, \cdot\right], \lhook\right)$ is a post-Lie algebra. The bracket $\llbracket-,-\rrbracket$ defines the Lie algebra $\left(\mathfrak{g}, \llbracket-,-\rrbracket\right)$ of the affine feedback group $\left(G, \star\right)$, where
\begin{align}\label{eq:derived-bracket}
	\llbracket \eta \mathrm{e}_i, \zeta \mathrm{e}_j \rrbracket
	&= \left(\eta \mathrm{e}_i \lhook \zeta \mathrm{e}_j\right)
		- \left(\zeta \mathrm{e}_j \lhook \eta \mathrm{e}_i\right)
		+ \left[\eta \mathrm{e}_i, \zeta \mathrm{e}_j\right]
\end{align}
for all $\eta,\zeta \in X^{+} \cup \{\emptyset\}$ and $i,j = 1,2$.
\endth
\medskip

The differential of the Grossman--Larson group product $\star$ is denoted by
\begin{align*}
	\bullet : \mathfrak{g} \otimes \mathfrak{g} \longrightarrow \mathfrak{g}.
\end{align*}
For all $t,s \in \re$ and $\mathrm{u},\mathrm{v} \in \mathfrak{g}$ it follows that
\begin{align*}
\exp\left(t\mathrm{u}\right) \star \exp\left(s\mathrm{v}\right) &=	\left(\mathbf{e} + t\mathrm{u} + \mathfrak{o}(t)\right) \star \left(\mathbf{e} + s\mathrm{v} + \mathfrak{o}(s)\right) \\
&= \mathbf{e} + t\mathrm{u} + s\mathrm{v} + st\left(\mathrm{u}\bullet\mathrm{v}\right) + \mathfrak{o}(t) + \mathfrak{o}(s) + \mathfrak{o}(st). 
\end{align*}
\medskip

\begth\label{th:GL-prod-lin-ind}\cite{Foissy_18}
If $\eta, \zeta \in X^{+} \cup \{\emptyset\}$, then for $i,j = 1,2$
\begin{align}\label{eq:GL-grp-prod-lin}
	\eta \mathrm{e}_i \bullet \zeta \mathrm{e}_j
	&= \left(\eta \mathrm{e}_i \lhook \zeta \mathrm{e}_j\right)
		+ \left(\eta \shuffle \zeta\right)\mathrm{e}_j \delta_{i,1}.
\end{align}
\endth

\begpr Let $t,s \in \re$. Then
\begin{align*}
\exp\left(t\eta\mathrm{e}_i\right) \star \exp\left(s\zeta\mathrm{e}_j\right)	&= \left(\mathbf{e} + t\eta \mathrm{e}_i + \mathfrak{o}(t)\right) \star \left(\mathbf{e} + s \zeta e_j + \mathfrak{o}(s)\right) \\
	&\stackrel{\eqref{eqn:grp-prod-affine-cord}}{=} \left[\left(\mathbf{e}
		+ t\eta \mathrm{e}_i +\mathfrak{o}(t)\right)\triangleleft \left(\mathbf{e}
		+ s\zeta \mathrm{e}_j + \mathfrak{o}(s)\right)\right] \cdot \left(\mathbf{e} + s\zeta \mathrm{e}_j + \mathfrak{o}(s)\right) \\
&\stackrel{\eqref{eqn:mixcomp-lin}}{=} \left[\mathbf{e} + t\eta \mathrm{e}_i
		+ st \left(\eta \mathrm{e}_i \lhook \zeta \mathrm{e}_j\right)
		+ \mathfrak{o}(t) + \mathfrak{o}(st)\right] \cdot \left(\mathbf{e} + s\zeta \mathrm{e}_j + \mathfrak{o}(s)\right) \\
	&\stackrel{\eqref{eqn:grp-prod-odot}}{=} \mathbf{e} + t\eta \mathrm{e}_i + s \zeta \mathrm{e}_j
	+ st \left[\left(\eta \mathrm{e}_i \lhook \zeta \mathrm{e}_j\right)
	+ \left(\eta \shuffle \zeta\right)\mathrm{e}_j \delta_{i,1}\right] \\
& \quad \quad \quad	+ \mathfrak{o}(t) + \mathfrak{o}(s) + \mathfrak{o}(st).
\end{align*}
\endpr
\medskip

As an independent check of Theorem~\ref{th:derived-Lie-affine-feedb}, observe that for all $\eta \mathrm{e}_i, \zeta \mathrm{e}_j \in \mathfrak{g}$:
\begin{align*}
	\left(\eta \mathrm{e}_i \bullet \zeta \mathrm{e}_j \right) - \left(\zeta \mathrm{e}_j \bullet \eta \mathrm{e}_i \right)
	&\stackrel{\eqref{eq:GL-grp-prod-lin}}{=} \left(\eta \mathrm{e}_i \lhook \zeta \mathrm{e}_j\right)
			- \left(\zeta \mathrm{e}_j \lhook \eta \mathrm{e}_i\right)
		+ \lbrace \left(\eta \shuffle \zeta\right)\mathrm{e}_j \delta_{i,1}
		- \left(\eta \shuffle \zeta\right)\mathrm{e}_i \delta_{j,1}\rbrace\\
	&\stackrel{\eqref{eq:Lie-bracket-def}}{=} \left(\eta \mathrm{e}_i \lhook \zeta \mathrm{e}_j\right)
	- \left(\zeta \mathrm{e}_j \lhook \eta \mathrm{e}_i\right)
	+ \left[\eta \mathrm{e}_i, \zeta \mathrm{e}_j\right] \\
	&\stackrel{\eqref{eq:derived-bracket}}{=} \llbracket \eta \mathrm{e}_i, \zeta \mathrm{e}_j\rrbracket.
\end{align*}

\begre{
It is shown in \cite{Foissy_18} and \cite{Gray-KEF_SIAM2017} that $\left(\mathfrak{g}, \bullet\right)$ is a pre-Lie algebra. Thus, the Lie algebra of the affine feedback group can be derived from both a post-Lie structure and a pre-Lie structure. The general question of precisely when a given Lie algebra $\mathfrak{g}$ can be derived from both pre- and post-Lie structures is an open problem.
}


\section{$\mathcal{H}$ and the Cointeraction}
\label{sec:hopf-alg}

Let $\mathcal{H}$ be the Hopf algebra of coordinate functions on $G$. The bialgebra of coordinate functions on $G$ with the coproduct dualizing the group product $\star$ was developed in \cite{Gray-KEF_SIAM2017}.\footnote{The coordinate functions are denoted by $a_{\eta}$ and $b_{\eta}$ instead of $\eta\epsilon_1, \eta \epsilon_2$ as in \cite{Foissy_18} and the current document.} The contribution of the present work is to describe another coproduct $\delta$ in $\mathcal{H}$ dualizing the group product $\cdot$. Further dualizing the action map, $\triangleleft$, one obtains the coaction map leading to a cointeraction between bialgebras $\left(\mathcal{H},\delta\right)$ and $\left(\mathcal{H},\Delta\right)$. For more details on cointeractig bialgebras, refer~\cite{Foissy_22}.
\medskip

The algebra structure of $\mathcal{H}$ is described briefly next and the coproduct $\delta$ is defined.

\begen[(i)]
\item Let $\allpoly{}^2$ be the vector space spanned by $x_{i_1}x_{i_2}\cdots x_{i_k}\epsilon_j$ for all $k \geq 0$ and $j=1,2$. The empty word ($k=0$) is denoted by $1$. For $\epsilon_j$ it follows that $\epsilon_j(\mathrm{e}_i) = \delta_{i,j}$

\item Let $\mathcal{S}\left(\allpoly{}^2\right)$ be the free symmetric algebra over $\allpoly{}^2$ with unit $\mbf{1}$.

\item For $\mathbf{c} = c_1\mathrm{e}_1 + c_2\mathrm{e}_2 \in \allseries{}^2$ and any $\eta \in \allwords$:
$1 \epsilon_i (\mathbf{c}) = c_i(1)$, $\eta \epsilon_i (\mathbf{c}) = c_i(\eta)$, $i = 1,2$.

\item As an algebra, $\mathcal{H} \cong \mathcal{S}\left(\allpoly{}^2\right)/\left(1\epsilon_1 - \mbf{1}\right)$. The symmetric product is denoted by $\HAprod$.

\item Define the coproduct $\delta$  dualizing the group product $\cdot$ on the coordinate functions to be
\begin{align*}
	\delta\left(\eta \epsilon_i\right)\left(\mathbf{c} \otimes \mathbf{d}\right)
	&= \left(\mathbf{c}\cdot\mathbf{d}\right)_{i}\left(\eta\right), \quad i=1,2,
\end{align*}
$\eta \in \allwords$ and $\mathbf{c},\mathbf{d} \in G$. The counit $\varepsilon$ is defined as $\varepsilon\left(\mbf{1}\right) = 1$ and zero otherwise.

\item For all $\eta \in \allwords$:
\begin{align*}
\begin{aligned}
	\delta \left(\eta \epsilon_1\right)
	&= \Delta_{\shuffle}\left(\eta\right)\left(\epsilon_1 \otimes \epsilon_1\right)\\
	\delta \left(\eta\epsilon_2\right)
	&= \left(\Delta_{\shuffle}\left(\eta\right)\left(\epsilon_1 \otimes \epsilon_2\right)\right) + \eta\epsilon_2 \otimes \mbf{1}.
\end{aligned}
\end{align*}

\item Let the degree of an element $\eta\epsilon_i$ be denoted by $\deg\left(\eta\epsilon_i\right)$. For all $x_{i_1}x_{i_2}\cdots x_{i_k} \in \allwords$:
\begin{align*}
	\deg\left(x_{i_1}x_{i_2}\cdots x_{i_k}\epsilon_1\right)
	&= 2\left[\sum_{j = 1}^k \delta_{i_j,0}\right] + \left[\sum_{j=1}^k \delta_{i_j,1}\right] \nonumber\\
	\deg\left(x_{i_1}x_{i_2}\cdots x_{i_k}\epsilon_2\right)
	& = 2\left[\sum_{j = 1}^k \delta_{i_j,0}\right] + \left[\sum_{j=1}^k \delta_{i_j,1}\right] + 1.
\end{align*}

\item Let $X_k = \{\eta\epsilon_i : \eta \in \allwords, \; \deg\left(\eta\epsilon_i\right) = k\}$ and $\mathcal{V}_k = \text{Vect}\left(X_k\right)$. The graded vector space $\mathcal{V} = \bigoplus_{n \geq 0} \mathcal{V}_n$ induces the grading on: $\mathcal{H} = \bigoplus_{n \geq 0} \mathcal{H}_n.$
Note that $\mathcal{H}_0 \cong \re$.

\item It is straightforward to check that for all $n \geq 0$: $\delta\left(\mathcal{V}_n\right) \subseteq \left(\mathcal{V} \otimes \mathcal{V}\right)_n.$ Therefore, $\left(\mathcal{H},\HAprod,\mbf{1},\delta,\varepsilon\right)$ is a commutative graded connected bialgebra.
\enden


\subsection{The Coaction Map $\rho$ and Coproduct $\Delta$}

Define a unital algebra morphism
\begin{align*}
	\rho : \mathcal{H} &\longrightarrow \mathcal{H} \otimes \mathcal{H} \nonumber \\
	\eta\epsilon_i &\longmapsto \rho\left(\eta\epsilon_i\right)\left(\mathbf{c}\otimes\mathbf{d}\right)
	= \left(\mathbf{c} \triangleleft \mathbf{d}\right)_{i}\left(\eta\right), \quad \;i=1,2,\;\; \mathbf{c},\mathbf{d} \in G.
\end{align*}
Define the linear maps
\begin{align*}
	\theta_j : \mathcal{V} &\longrightarrow \mathcal{V} \\
 	\eta \epsilon_i &\longmapsto x_j\eta \epsilon_i, \quad j = 0,1,\;\; i=1,2.
\end{align*}
The computations for $\rho$ on $\mathcal{V}$ are described in \cite{Gray-KEF_SIAM2017}\footnote{In \cite{Gray-KEF_SIAM2017}, the notation is $\tilde{\Delta}$.} and Foissy \cite{Foissy_18}\footnote{In \cite{Foissy_18}, this was called the reduced coproduct with the notation being $\bar{\Delta}_{\ast}$.}. Specifically, for any word $\eta \in \allwords$ and $i =1,2$:
\begin{align*}
\begin{aligned}
	\rho\left(1\epsilon_2\right)
	&= 1\epsilon_2 \otimes \mbf{1}\\
	\rho \circ \theta_0\left(\eta\epsilon_i\right)
	&= \left(\theta_0 \otimes \id_{\mathcal{H}}\right) \circ \rho\left(\eta \epsilon_i\right)
	+ \left(\theta_1 \otimes \HAprod\right) \circ \left(\rho \otimes \id_{\mathcal{H}}\right)
	\Delta_{\shuffle}\left(\eta\right)\left(\epsilon_i \otimes \epsilon_2\right)  \\
	\rho \circ \theta_1\left(\eta\epsilon_i\right)
	&= \left(\theta_1 \otimes \HAprod \right) \circ \left(\rho \otimes \id_{\mathcal{H}}\right)
	\Delta_{\shuffle}\left(\eta\right)\left(\epsilon_i \otimes \epsilon_1\right).
\end{aligned}
\end{align*}

\begth\cite{Gray-KEF_SIAM2017}\footnote{See Lemma $4.3$ of \cite{Gray-KEF_SIAM2017}.}
The morphism $\rho$ respects the grading on $\mathcal{H}$ viz.~$\rho\left(\mathcal{V}_n\right) \subseteq \left(\mathcal{V} \otimes \mathcal{H}\right)_n$. Hence, $\rho : \mathcal{H} \longrightarrow \mathcal{H} \otimes \mathcal{H}$ is a graded unital algebra morphism.
\endth

Define the coproduct $\Delta : \mathcal{H} \longrightarrow \mathcal{H} \otimes \mathcal{H}$ dualizing the Grossman--Larson group product $\star$
\begin{align*}
	\Delta\left(\eta\epsilon_i\right)\left(\mathbf{c} \otimes \mathbf{d}\right)
	= \left(\mathbf{c}\star\mathbf{d}\right)_i \left(\eta\right), \quad i = 1,2, \;\; \eta \in \allwords, \;\; \mathbf{c}, \mathbf{d} \in G.
\end{align*}

\begth\cite{Gray-KEF_SIAM2017}\footnote{See Theorem~$4.6$ of \cite{Gray-KEF_SIAM2017}.}
$\left(\mathcal{H},\HAprod,\mbf{1},\Delta,\epsilon\right)$ is a connected graded commutative bialgebra.
\endth

The computation of the coproduct $\Delta$ can be expressed using $\delta$ and $\rho$ dualizing \rref{eqn:grp-prod-affine-cord}:
\begin{align}
\label{eq:feedb-rel}
	\Delta &= \left(\id_{\mathcal{H}} \otimes \HAprod\right) \circ \left(\rho \otimes \id_{\mathcal{H}}\right) \circ \delta.
\end{align}
Dualizing Theorem~\ref{th:affine_feedb} on the coordinate functions of $G$, one sees that $\left(\mathcal{H},\HAprod,\delta,\varepsilon \right)$ is a Hopf algebra in the category of right $\left(\mathcal{H}, \HAprod,\Delta, \epsilon\right)$-comodule via the coaction map $\rho$.
\medskip

\begre{~The cointeraction in $\mathcal{H}$ is novel to the current work and was not obsevred either in \cite{Gray-KEF_SIAM2017} or \cite{Foissy_18}. For in \cite{Foissy_18}, the coaction $\rho$ was termed as ``reduced coproduct" and both the mentioned works did not notice the second coalgebra endowed by the coproduct $\delta$.}


\section{Concluding Remarks}
\label{sec:concluding-rem}

The following arrow diagram summarized the algebraic structures described in this paper.

\begin{center}
\label{eq:arrow-diag}
\begin{tikzcd}[sep = 14ex]
\left(\mathcal{H},\delta\right)\xleftharpoondown[]{\rho}\left(\mathcal{H},\Delta\right)  \arrow[ddr, bend right = 10, "(6)"] \arrow[dd,"(3)"] &  &  \\
& \left(G,\cdot\right)\xleftharpoonup[]{\triangleleft}\left(G,\star\right) \arrow[ul, "(1\&2)"] \arrow[d, "(5)"]  \arrow[dr, "(8)"] &  \\
\left(\mathfrak{g}_{\mathcal{SISO}},\left[-,-\right],\lhook\right) \arrow[r, "(4)"] &  \left(\mathfrak{g},\left[-,-\right],\lhook\right) \arrow[r, "(7)"] &  \left(\mathcal{U}\left(\mathfrak{g}\right), \oast, \Delta_{\shuffle}\right)
\end{tikzcd}
\end{center}

\begen
\item  $\left(G,\cdot\right)\xleftharpoonup[]{\triangleleft}\left(G,\star\right)$: The affine feedback group $\left(G,\star\right)$ that appears in the study of affine feedback interconnection of Chen--Fliess series is shown to be the Grossman--Larson group of the post-group $\left(G,\cdot,\triangleleft\right)$, which is the central object in this paper. It acts on the group $\left(G, \cdot\right)$.
\medskip

\item $\left(\mathcal{H},\delta\right)\xleftharpoondown[]{\rho}\left(\mathcal{H},\Delta\right)$: The coordinate functions on the post-group $\left(G,\cdot,\triangleleft\right)$ form a cointeracting bialgebra structure $\left(\mathcal{H}, \HAprod,\mbf{1},\Delta,\delta,\rho,\Delta,\epsilon\right)$, where the coproducts $\Delta$ and $\delta$ dualize the group products $\star$ and $\cdot$, respectively, while the coaction $\rho$ dualizes $\triangleleft$. $\left(\mathcal{H},\HAprod, \delta\right)$ is a bialgebra in the category of $\left(\mathcal{H},\HAprod, \delta\right)$-comodules via the coaction map $\rho$. 
The coproducts $\Delta$, $\delta$ and the coaction map $\rho$ are related by the relation~\eqref{eq:feedb-rel}, which dualizes the Guin--Oudom construction of the Grossman--Larson group product.

\item $\mathfrak{g}_{\scriptscriptstyle{\mathcal{SISO}}}$: In \cite{Foissy_18}, Foissy linearized the coproduct $\Delta$ and coaction $\rho$. The post-Lie algebra $\left(\mathfrak{g}_{\scriptscriptstyle{\mathcal{SISO}}},\left[-,-\right],\lhook\right)$ followed from taking the graded dual of the linearized coproducts and coaction maps.

\item $\mathfrak{g}$: The filtered completion of the post-Lie algebra $\mathfrak{g}_{\scriptscriptstyle{\mathcal{SISO}}}$ is the post-Lie algebra $\left(\mathfrak{g},\left[-,-\right],\lhook\right)$.

\item The completed post-Lie algebra $\mathfrak{g}$ can be obtained from the Hopf algebra $\mathcal{H}$ by linearizing the coproduct $\delta$ and coaction map $\rho$ as in \cite{Foissy_18} and then taking their linear duals (instead of the graded dual). The derived Lie algebra $\left(\mathfrak{g}, \llbracket-,- \rrbracket\right)$ is the Lie algebra corresponding to the Grossman--Larson group $\left(G, \star\right)$.

\item The post-Lie algebra $\left(\mathfrak{g}, \left[-,-\right], \lhook\right)$ can also be derived by differentiating the group product and  action in $\left(G,\cdot,\triangleleft\right)$.

\item $\left(\mathcal{U}(\mathfrak{g}), \oast, \Delta_{\shuffle}\right)$: Let $\mathcal{U}(\mathfrak{g})$ denote the universal enveloping algebra of the Lie algebra $\left(\mathfrak{g},\left[-,-\right]\right)$. Extending the post-Lie product $\lhook$ to $\mathcal{U}(\mathfrak{g})$ \cite{MuntheKaas2015postLie} results in the bialgebra $\left(\mathcal{U}(\mathfrak{g}), \oast, \Delta_{\shuffle}\right)$, which is isomorphic (as a graded Hopf algebra) to the universal enveloping algebra of $\left(\mathfrak{g}, \llbracket-,- \rrbracket \right)$ \cite{MuntheKaas2015postLie}.

\item The group elements in $\left(\mathcal{U}(\mathfrak{g}), \oast, \Delta_{\shuffle}\right)$ are isomorphic to the elements of the Grossman--Larson group $\left(G, \star\right)$.
\enden


\end{document}